\documentclass{gtart}

\def\ifplaintex{\expandafter\ifx\csname documentclass\endcsname\relax}

\def\gtp{{\mathsurround=0pt\it $\cal G\mskip-2mu$eometry \&\ 
$\cal T\!\!$opology $\cal P\!$ublications}}  

\def\recd{{\small Received:\qua\receiveddate\ifx\reviseddate\relax
\else\qquad Revised:\qua\reviseddate\fi\par}} 

\def\lognumber#1{\def\thelognumber{#1}}
\def\volumenumber#1{\def\thevolumenumber{#1}}
\def\volumeyear#1{\def\thevolumeyear{#1}}
\def\papernumber#1{\def\thepapernumber{#1}}
\def\pagenumbers#1#2{\def\startpage{#1}\def\finishpage{#2}}
\def\published#1{\def\publishdate{#1}}

\def\received#1{\def\receiveddate{#1}}
\def\revised#1{\def\reviseddate{#1}}
\def\accepted#1{\def\accepteddate{#1}}
\def\asciititle#1{\def\theasciititle{#1}}

\let\\\par\let\thelognumber\relax\let\thevolumenumber\relax
\let\thepapernumber\relax\let\thevolumeyear\relax\let\startpage\relax
\let\finishpage\relax\let\publishdate\relax\let\receiveddate\relax
\let\reviseddate\relax\let\accepteddate\relax\let\theasciititle\relax
\let\theasciiauthors\relax
\let\theasciiabstract\relax

\let\theasciiemail\relax

\ifplaintex
\font\logobig=cmssbx10 scaled 3836
\font\logomed=cmssbx10 scaled 2557
\else
\font\logobig=cmssbx10 scaled 4200
\font\logomed=cmssbx10 scaled 2800
\fi

\long\def\makeagttitle{   
\count0=\startpage
\agt\hfill      
\hbox to 45truept{\vbox to 0pt{\vglue -13truept{\logomed A\kern -.37em{\logobig 
T}\kern -.38em G}\vss}\hss}
\break
{\small Volume \thevolumenumber\ (\thevolumeyear)
\startpage--\finishpage\nl
Published: \publishdate}

\vglue .25truein

{\parskip=0pt\leftskip 0pt plus
1fil\def\\{\par\smallskip}{\Large\bf\thetitle}\par\medskip} \vglue
0.05truein

{\parskip=0pt\leftskip 0pt plus 1fil\def\\{\par}{\sc\theauthors}
\par\medskip}%
 
\vglue 0.03truein 

{\small\leftskip 25truept\rightskip 25truept{\bf Abstract}\stdspace\theabstract

{\bf AMS Classification}\stdspace\theprimaryclass
\ifx\thesecondaryclass\relax\else; \thesecondaryclass\fi\par
{\bf Keywords}\stdspace \thekeywords\par}\vglue 7truept

}   

\ifplaintex
\font\phead=cmsl9 scaled 950
\font\pnum=cmbx10 scaled 913
\font\pfoot=cmsl9 scaled 950
\headline{\vbox to 0pt{\vskip -4.5mm\line{\small\phead\ifnum
\count0=\startpage ISSN 1472-2739 (on-line) 1472-2747 (printed)
\hfill {\pnum\folio}\else\ifodd\count0\def\\{ }%
\ifx\theshorttitle\relax\thetitle\else\theshorttitle\fi\hfill{\pnum\folio}
\else\def\\{ and }{\pnum\folio}\hfill\ifx\theshortauthors\relax\theauthors
\else\theshortauthors\fi\fi\fi}\vss}}
\footline{\vbox to 0pt{\vglue 0mm\line{\small\pfoot\ifnum\count0=\startpage
\copyright\ \gtp\hfill\else
\agt, Volume \thevolumenumber\ (\thevolumeyear)\hfill\fi}\vss}}
\else
\font=cmsl9 scaled 1050
\font\lnum=cmbx10 
\font=cmsl9 scaled 1050
\makeatletter
\def\@oddhead{{\small\lhead\ifnum\count0=\startpage ISSN 1472-2739 
(on-line) 1472-2747 (printed)\hfill {\lnum\number\count0}\else\ifodd\count0
\def\\{ }\ifx\theshorttitle\relax \thetitle \else\theshorttitle\fi\hfill
{\lnum\number\count0}\else\def\\{ and }{\lnum\number\count0}
\hfill\ifx\theshortauthors\relax 
\theauthors\else\theshortauthors\fi\fi\fi}}\def\@evenhead{\@oddhead}
\def\@oddfoot{\small\lfoot\ifnum\count0=\startpage\copyright\ \gtp\hfill\else
\agt, Volume \thevolumenumber\ (\thevolumeyear)\hfill\fi}
\def\@evenfoot{\@oddfoot}
\makeatother
\fi
\let\maketitlepage\makeagttitle
\let\makeshorttitle\maketitlepage
\let\maketitle\maketitlepage

\newwrite\gtoutfile
\long\gdef\makeheadfile{  
{\def\\{, }\def\s{ }
\immediate\openout\gtoutfile head.xxx
\immediate\write\gtoutfile{To: math@arxiv.org}
\immediate\write\gtoutfile{Subject: put OR rep NNNNN:ppppp}
\immediate\write\gtoutfile{--text follows this line--}
\immediate\write\gtoutfile{Proxy-for: \ifx\theasciiauthors\relax
\theauthors\else\theasciiauthors\fi\s<\ifx\theasciiemail\relax\theemail\else\theasciiemail\fi>}
\immediate\write\gtoutfile{\noexpand\\}
\immediate\write\gtoutfile{Authors: \ifx\theasciiauthors\relax
\theauthors\else\theasciiauthors\fi}
{\def\\{ }\immediate\write\gtoutfile{Title: \ifx\theasciititle\relax
\thetitle\else\theasciititle\fi}}
\immediate\write\gtoutfile{Subj-class: GT or SG, GR etc}
\immediate\write\gtoutfile{MSC-class: \theprimaryclass\ifx\thesecondaryclass\relax\else, \thesecondaryclass\fi}
\immediate\write\gtoutfile{Journal-ref: Algebr. Geom. Topol. \thevolumenumber\s
(\thevolumeyear) \startpage-\finishpage}
\immediate\write\gtoutfile{Comments: Published by Algebraic and
Geometric Topology at}
\immediate\write\gtoutfile{\s\s\s  http://www.maths.warwick.ac.uk/agt/AGTVol\thevolumenumber/agt-\thevolumenumber-\thepapernumber.abs.html}
\immediate\write\gtoutfile{\noexpand\\}
\immediate\write\gtoutfile{}
\ifx\theasciiabstract\relax
\immediate\write\gtoutfile{\theabstract}\else
\immediate\write\gtoutfile{\theasciiabstract}\fi
\immediate\write\gtoutfile{}
\immediate\write\gtoutfile{\noexpand\\}
\immediate\write\gtoutfile{}
\immediate\closeout\gtoutfile}}  

\def\maketitlepage{\makeagttitle\makeheadfile}
\let\makeshorttitle\maketitlepage
\let\maketitle\maketitlepage

\def\ifplaintex{\expandafter\ifx\csname documentclass\endcsname\relax}

\def\gtp{{\mathsurround=0pt\it $\cal G\mskip-2mu$eometry \&\ 
$\cal T\!\!$opology $\cal P\!$ublications}}  

\def\recd{{\small Received:\qua\receiveddate\ifx\reviseddate\relax
\else\qquad Revised:\qua\reviseddate\fi\par}} 

\def\lognumber#1{\def\thelognumber{#1}}
\def\volumenumber#1{\def\thevolumenumber{#1}}
\def\volumeyear#1{\def\thevolumeyear{#1}}
\def\papernumber#1{\def\thepapernumber{#1}}
\def\pagenumbers#1#2{\def\startpage{#1}\def\finishpage{#2}}
\def\published#1{\def\publishdate{#1}}

\def\received#1{\def\receiveddate{#1}}
\def\revised#1{\def\reviseddate{#1}}
\def\accepted#1{\def\accepteddate{#1}}
\def\asciititle#1{\def\theasciititle{#1}}

\let\\\par\let\thelognumber\relax\let\thevolumenumber\relax
\let\thepapernumber\relax\let\thevolumeyear\relax\let\startpage\relax
\let\finishpage\relax\let\publishdate\relax\let\receiveddate\relax
\let\reviseddate\relax\let\accepteddate\relax\let\theasciititle\relax
\let\theasciiauthors\relax
\let\theasciiabstract\relax

\let\theasciiemail\relax

\ifplaintex
\font\logobig=cmssbx10 scaled 3836
\font\logomed=cmssbx10 scaled 2557
\else
\font\logobig=cmssbx10 scaled 4200
\font\logomed=cmssbx10 scaled 2800
\fi

\long\def\makeagttitle{   
\count0=\startpage
\agt\hfill      
\hbox to 45truept{\vbox to 0pt{\vglue -13truept{\logomed A\kern -.37em{\logobig 
T}\kern -.38em G}\vss}\hss}
\break
{\small Volume \thevolumenumber\ (\thevolumeyear)
\startpage--\finishpage\nl
Published: \publishdate}

\vglue .25truein

{\parskip=0pt\leftskip 0pt plus
1fil\def\\{\par\smallskip}{\Large\bf\thetitle}\par\medskip} \vglue
0.05truein

{\parskip=0pt\leftskip 0pt plus 1fil\def\\{\par}{\sc\theauthors}
\par\medskip}%
 
\vglue 0.03truein 

{\small\leftskip 25truept\rightskip 25truept{\bf Abstract}\stdspace\theabstract

{\bf AMS Classification}\stdspace\theprimaryclass
\ifx\thesecondaryclass\relax\else; \thesecondaryclass\fi\par
{\bf Keywords}\stdspace \thekeywords\par}\vglue 7truept

}   

\ifplaintex
\font\phead=cmsl9 scaled 950
\font\pnum=cmbx10 scaled 913
\font\pfoot=cmsl9 scaled 950
\headline{\vbox to 0pt{\vskip -4.5mm\line{\small\phead\ifnum
\count0=\startpage ISSN 1472-2739 (on-line) 1472-2747 (printed)
\hfill {\pnum\folio}\else\ifodd\count0\def\\{ }%
\ifx\theshorttitle\relax\thetitle\else\theshorttitle\fi\hfill{\pnum\folio}
\else\def\\{ and }{\pnum\folio}\hfill\ifx\theshortauthors\relax\theauthors
\else\theshortauthors\fi\fi\fi}\vss}}
\footline{\vbox to 0pt{\vglue 0mm\line{\small\pfoot\ifnum\count0=\startpage
\copyright\ \gtp\hfill\else
\agt, Volume \thevolumenumber\ (\thevolumeyear)\hfill\fi}\vss}}
\else
\font=cmsl9 scaled 1050
\font\lnum=cmbx10 
\font=cmsl9 scaled 1050
\makeatletter
\def\@oddhead{{\small\lhead\ifnum\count0=\startpage ISSN 1472-2739 
(on-line) 1472-2747 (printed)\hfill {\lnum\number\count0}\else\ifodd\count0
\def\\{ }\ifx\theshorttitle\relax \thetitle \else\theshorttitle\fi\hfill
{\lnum\number\count0}\else\def\\{ and }{\lnum\number\count0}
\hfill\ifx\theshortauthors\relax 
\theauthors\else\theshortauthors\fi\fi\fi}}\def\@evenhead{\@oddhead}
\def\@oddfoot{\small\lfoot\ifnum\count0=\startpage\copyright\ \gtp\hfill\else
\agt, Volume \thevolumenumber\ (\thevolumeyear)\hfill\fi}
\def\@evenfoot{\@oddfoot}
\makeatother
\fi
\let\maketitlepage\makeagttitle
\let\makeshorttitle\maketitlepage
\let\maketitle\maketitlepage

\newwrite\gtoutfile
\long\gdef\makeheadfile{  
{\def\\{, }\def\s{ }
\immediate\openout\gtoutfile head.xxx
\immediate\write\gtoutfile{To: math@arxiv.org}
\immediate\write\gtoutfile{Subject: put OR rep NNNNN:ppppp}
\immediate\write\gtoutfile{--text follows this line--}
\immediate\write\gtoutfile{Proxy-for: \ifx\theasciiauthors\relax
\theauthors\else\theasciiauthors\fi\s<\ifx\theasciiemail\relax\theemail\else\theasciiemail\fi>}
\immediate\write\gtoutfile{\noexpand\\}
\immediate\write\gtoutfile{Authors: \ifx\theasciiauthors\relax
\theauthors\else\theasciiauthors\fi}
{\def\\{ }\immediate\write\gtoutfile{Title: \ifx\theasciititle\relax
\thetitle\else\theasciititle\fi}}
\immediate\write\gtoutfile{Subj-class: GT or SG, GR etc}
\immediate\write\gtoutfile{MSC-class: \theprimaryclass\ifx\thesecondaryclass\relax\else, \thesecondaryclass\fi}
\immediate\write\gtoutfile{Journal-ref: Algebr. Geom. Topol. \thevolumenumber\s
(\thevolumeyear) \startpage-\finishpage}
\immediate\write\gtoutfile{Comments: Published by Algebraic and
Geometric Topology at}
\immediate\write\gtoutfile{\s\s\s  http://www.maths.warwick.ac.uk/agt/AGTVol\thevolumenumber/agt-\thevolumenumber-\thepapernumber.abs.html}
\immediate\write\gtoutfile{\noexpand\\}
\immediate\write\gtoutfile{}
\ifx\theasciiabstract\relax
\immediate\write\gtoutfile{\theabstract}\else
\immediate\write\gtoutfile{\theasciiabstract}\fi
\immediate\write\gtoutfile{}
\immediate\write\gtoutfile{\noexpand\\}
\immediate\write\gtoutfile{}
\immediate\closeout\gtoutfile}}  

\def\maketitlepage{\makeagttitle\makeheadfile}
\let\makeshorttitle\maketitlepage
\let\maketitle\maketitlepage

\def\ifplaintex{\expandafter\ifx\csname documentclass\endcsname\relax}

\def\gtp{{\mathsurround=0pt\it $\cal G\mskip-2mu$eometry \&\ 
$\cal T\!\!$opology $\cal P\!$ublications}}  

\def\recd{{\small Received:\qua\receiveddate\ifx\reviseddate\relax
\else\qquad Revised:\qua\reviseddate\fi\par}} 

\def\lognumber#1{\def\thelognumber{#1}}
\def\volumenumber#1{\def\thevolumenumber{#1}}
\def\volumeyear#1{\def\thevolumeyear{#1}}
\def\papernumber#1{\def\thepapernumber{#1}}
\def\pagenumbers#1#2{\def\startpage{#1}\def\finishpage{#2}}
\def\published#1{\def\publishdate{#1}}

\def\received#1{\def\receiveddate{#1}}
\def\revised#1{\def\reviseddate{#1}}
\def\accepted#1{\def\accepteddate{#1}}
\def\asciititle#1{\def\theasciititle{#1}}

\let\\\par\let\thelognumber\relax\let\thevolumenumber\relax
\let\thepapernumber\relax\let\thevolumeyear\relax\let\startpage\relax
\let\finishpage\relax\let\publishdate\relax\let\receiveddate\relax
\let\reviseddate\relax\let\accepteddate\relax\let\theasciititle\relax
\let\theasciiauthors\relax
\let\theasciiabstract\relax

\let\theasciiemail\relax

\ifplaintex
\font\logobig=cmssbx10 scaled 3836
\font\logomed=cmssbx10 scaled 2557
\else
\font\logobig=cmssbx10 scaled 4200
\font\logomed=cmssbx10 scaled 2800
\fi

\long\def\makeagttitle{   
\count0=\startpage
\agt\hfill      
\hbox to 45truept{\vbox to 0pt{\vglue -13truept{\logomed A\kern -.37em{\logobig 
T}\kern -.38em G}\vss}\hss}
\break
{\small Volume \thevolumenumber\ (\thevolumeyear)
\startpage--\finishpage\nl
Published: \publishdate}

\vglue .25truein

{\parskip=0pt\leftskip 0pt plus
1fil\def\\{\par\smallskip}{\Large\bf\thetitle}\par\medskip} \vglue
0.05truein

{\parskip=0pt\leftskip 0pt plus 1fil\def\\{\par}{\sc\theauthors}
\par\medskip}%
 
\vglue 0.03truein 

{\small\leftskip 25truept\rightskip 25truept{\bf Abstract}\stdspace\theabstract

{\bf AMS Classification}\stdspace\theprimaryclass
\ifx\thesecondaryclass\relax\else; \thesecondaryclass\fi\par
{\bf Keywords}\stdspace \thekeywords\par}\vglue 7truept

}   

\ifplaintex
\font\phead=cmsl9 scaled 950
\font\pnum=cmbx10 scaled 913
\font\pfoot=cmsl9 scaled 950
\headline{\vbox to 0pt{\vskip -4.5mm\line{\small\phead\ifnum
\count0=\startpage ISSN 1472-2739 (on-line) 1472-2747 (printed)
\hfill {\pnum\folio}\else\ifodd\count0\def\\{ }%
\ifx\theshorttitle\relax\thetitle\else\theshorttitle\fi\hfill{\pnum\folio}
\else\def\\{ and }{\pnum\folio}\hfill\ifx\theshortauthors\relax\theauthors
\else\theshortauthors\fi\fi\fi}\vss}}
\footline{\vbox to 0pt{\vglue 0mm\line{\small\pfoot\ifnum\count0=\startpage
\copyright\ \gtp\hfill\else
\agt, Volume \thevolumenumber\ (\thevolumeyear)\hfill\fi}\vss}}
\else
\font=cmsl9 scaled 1050
\font\lnum=cmbx10 
\font=cmsl9 scaled 1050
\makeatletter
\def\@oddhead{{\small\lhead\ifnum\count0=\startpage ISSN 1472-2739 
(on-line) 1472-2747 (printed)\hfill {\lnum\number\count0}\else\ifodd\count0
\def\\{ }\ifx\theshorttitle\relax \thetitle \else\theshorttitle\fi\hfill
{\lnum\number\count0}\else\def\\{ and }{\lnum\number\count0}
\hfill\ifx\theshortauthors\relax 
\theauthors\else\theshortauthors\fi\fi\fi}}\def\@evenhead{\@oddhead}
\def\@oddfoot{\small\lfoot\ifnum\count0=\startpage\copyright\ \gtp\hfill\else
\agt, Volume \thevolumenumber\ (\thevolumeyear)\hfill\fi}
\def\@evenfoot{\@oddfoot}
\makeatother
\fi
\let\maketitlepage\makeagttitle
\let\makeshorttitle\maketitlepage
\let\maketitle\maketitlepage

\newwrite\gtoutfile
\long\gdef\makeheadfile{  
{\def\\{, }\def\s{ }
\immediate\openout\gtoutfile head.xxx
\immediate\write\gtoutfile{To: math@arxiv.org}
\immediate\write\gtoutfile{Subject: put OR rep NNNNN:ppppp}
\immediate\write\gtoutfile{--text follows this line--}
\immediate\write\gtoutfile{Proxy-for: \ifx\theasciiauthors\relax
\theauthors\else\theasciiauthors\fi\s<\ifx\theasciiemail\relax\theemail\else\theasciiemail\fi>}
\immediate\write\gtoutfile{\noexpand\\}
\immediate\write\gtoutfile{Authors: \ifx\theasciiauthors\relax
\theauthors\else\theasciiauthors\fi}
{\def\\{ }\immediate\write\gtoutfile{Title: \ifx\theasciititle\relax
\thetitle\else\theasciititle\fi}}
\immediate\write\gtoutfile{Subj-class: GT or SG, GR etc}
\immediate\write\gtoutfile{MSC-class: \theprimaryclass\ifx\thesecondaryclass\relax\else, \thesecondaryclass\fi}
\immediate\write\gtoutfile{Journal-ref: Algebr. Geom. Topol. \thevolumenumber\s
(\thevolumeyear) \startpage-\finishpage}
\immediate\write\gtoutfile{Comments: Published by Algebraic and
Geometric Topology at}
\immediate\write\gtoutfile{\s\s\s  http://www.maths.warwick.ac.uk/agt/AGTVol\thevolumenumber/agt-\thevolumenumber-\thepapernumber.abs.html}
\immediate\write\gtoutfile{\noexpand\\}
\immediate\write\gtoutfile{}
\ifx\theasciiabstract\relax
\immediate\write\gtoutfile{\theabstract}\else
\immediate\write\gtoutfile{\theasciiabstract}\fi
\immediate\write\gtoutfile{}
\immediate\write\gtoutfile{\noexpand\\}
\immediate\write\gtoutfile{}
\immediate\closeout\gtoutfile}}  

\def\maketitlepage{\makeagttitle\makeheadfile}
\let\makeshorttitle\maketitlepage
\let\maketitle\maketitlepage

\def\ifplaintex{\expandafter\ifx\csname documentclass\endcsname\relax}

\def\gtp{{\mathsurround=0pt\it $\cal G\mskip-2mu$eometry \&\ 
$\cal T\!\!$opology $\cal P\!$ublications}}  

\def\recd{{\small Received:\qua\receiveddate\ifx\reviseddate\relax
\else\qquad Revised:\qua\reviseddate\fi\par}} 

\def\lognumber#1{\def\thelognumber{#1}}
\def\volumenumber#1{\def\thevolumenumber{#1}}
\def\volumeyear#1{\def\thevolumeyear{#1}}
\def\papernumber#1{\def\thepapernumber{#1}}
\def\pagenumbers#1#2{\def\startpage{#1}\def\finishpage{#2}}
\def\published#1{\def\publishdate{#1}}

\def\received#1{\def\receiveddate{#1}}
\def\revised#1{\def\reviseddate{#1}}
\def\accepted#1{\def\accepteddate{#1}}
\def\asciititle#1{\def\theasciititle{#1}}

\let\\\par\let\thelognumber\relax\let\thevolumenumber\relax
\let\thepapernumber\relax\let\thevolumeyear\relax\let\startpage\relax
\let\finishpage\relax\let\publishdate\relax\let\receiveddate\relax
\let\reviseddate\relax\let\accepteddate\relax\let\theasciititle\relax
\let\theasciiauthors\relax
\let\theasciiabstract\relax

\let\theasciiemail\relax

\ifplaintex
\font\logobig=cmssbx10 scaled 3836
\font\logomed=cmssbx10 scaled 2557
\else
\font\logobig=cmssbx10 scaled 4200
\font\logomed=cmssbx10 scaled 2800
\fi

\long\def\makeagttitle{   
\count0=\startpage
\agt\hfill      
\hbox to 45truept{\vbox to 0pt{\vglue -13truept{\logomed A\kern -.37em{\logobig 
T}\kern -.38em G}\vss}\hss}
\break
{\small Volume \thevolumenumber\ (\thevolumeyear)
\startpage--\finishpage\nl
Published: \publishdate}

\vglue .25truein

{\parskip=0pt\leftskip 0pt plus
1fil\def\\{\par\smallskip}{\Large\bf\thetitle}\par\medskip} \vglue
0.05truein

{\parskip=0pt\leftskip 0pt plus 1fil\def\\{\par}{\sc\theauthors}
\par\medskip}%
 
\vglue 0.03truein 

{\small\leftskip 25truept\rightskip 25truept{\bf Abstract}\stdspace\theabstract

{\bf AMS Classification}\stdspace\theprimaryclass
\ifx\thesecondaryclass\relax\else; \thesecondaryclass\fi\par
{\bf Keywords}\stdspace \thekeywords\par}\vglue 7truept

}   

\ifplaintex
\font\phead=cmsl9 scaled 950
\font\pnum=cmbx10 scaled 913
\font\pfoot=cmsl9 scaled 950
\headline{\vbox to 0pt{\vskip -4.5mm\line{\small\phead\ifnum
\count0=\startpage ISSN 1472-2739 (on-line) 1472-2747 (printed)
\hfill {\pnum\folio}\else\ifodd\count0\def\\{ }%
\ifx\theshorttitle\relax\thetitle\else\theshorttitle\fi\hfill{\pnum\folio}
\else\def\\{ and }{\pnum\folio}\hfill\ifx\theshortauthors\relax\theauthors
\else\theshortauthors\fi\fi\fi}\vss}}
\footline{\vbox to 0pt{\vglue 0mm\line{\small\pfoot\ifnum\count0=\startpage
\copyright\ \gtp\hfill\else
\agt, Volume \thevolumenumber\ (\thevolumeyear)\hfill\fi}\vss}}
\else
\font=cmsl9 scaled 1050
\font\lnum=cmbx10 
\font=cmsl9 scaled 1050
\makeatletter
\def\@oddhead{{\small\lhead\ifnum\count0=\startpage ISSN 1472-2739 
(on-line) 1472-2747 (printed)\hfill {\lnum\number\count0}\else\ifodd\count0
\def\\{ }\ifx\theshorttitle\relax \thetitle \else\theshorttitle\fi\hfill
{\lnum\number\count0}\else\def\\{ and }{\lnum\number\count0}
\hfill\ifx\theshortauthors\relax 
\theauthors\else\theshortauthors\fi\fi\fi}}\def\@evenhead{\@oddhead}
\def\@oddfoot{\small\lfoot\ifnum\count0=\startpage\copyright\ \gtp\hfill\else
\agt, Volume \thevolumenumber\ (\thevolumeyear)\hfill\fi}
\def\@evenfoot{\@oddfoot}
\makeatother
\fi
\let\maketitlepage\makeagttitle
\let\makeshorttitle\maketitlepage
\let\maketitle\maketitlepage

\newwrite\gtoutfile
\long\gdef\makeheadfile{  
{\def\\{, }\def\s{ }
\immediate\openout\gtoutfile head.xxx
\immediate\write\gtoutfile{To: math@arxiv.org}
\immediate\write\gtoutfile{Subject: put OR rep NNNNN:ppppp}
\immediate\write\gtoutfile{--text follows this line--}
\immediate\write\gtoutfile{Proxy-for: \ifx\theasciiauthors\relax
\theauthors\else\theasciiauthors\fi\s<\ifx\theasciiemail\relax\theemail\else\theasciiemail\fi>}
\immediate\write\gtoutfile{\noexpand\\}
\immediate\write\gtoutfile{Authors: \ifx\theasciiauthors\relax
\theauthors\else\theasciiauthors\fi}
{\def\\{ }\immediate\write\gtoutfile{Title: \ifx\theasciititle\relax
\thetitle\else\theasciititle\fi}}
\immediate\write\gtoutfile{Subj-class: GT or SG, GR etc}
\immediate\write\gtoutfile{MSC-class: \theprimaryclass\ifx\thesecondaryclass\relax\else, \thesecondaryclass\fi}
\immediate\write\gtoutfile{Journal-ref: Algebr. Geom. Topol. \thevolumenumber\s
(\thevolumeyear) \startpage-\finishpage}
\immediate\write\gtoutfile{Comments: Published by Algebraic and
Geometric Topology at}
\immediate\write\gtoutfile{\s\s\s  http://www.maths.warwick.ac.uk/agt/AGTVol\thevolumenumber/agt-\thevolumenumber-\thepapernumber.abs.html}
\immediate\write\gtoutfile{\noexpand\\}
\immediate\write\gtoutfile{}
\ifx\theasciiabstract\relax
\immediate\write\gtoutfile{\theabstract}\else
\immediate\write\gtoutfile{\theasciiabstract}\fi
\immediate\write\gtoutfile{}
\immediate\write\gtoutfile{\noexpand\\}
\immediate\write\gtoutfile{}
\immediate\closeout\gtoutfile}}  

\def\maketitlepage{\makeagttitle\makeheadfile}
\let\makeshorttitle\maketitlepage
\let\maketitle\maketitlepage

\lognumber{35}
\volumenumber{1}
\volumeyear{2001}
\papernumber{35}
\pagenumbers{709}{718}
\received{8 November 2001}
\revised{27 November 2001}
\accepted{28 November 2001}
\published{30 November 2001}

\usepackage{amssymb}
\usepackage{amscd}
\usepackage[leqno]{amsmath} 
\numberwithin{equation}{section} 

\newtheorem{Lemma}[subsection]{Lemma} 

\newtheorem{Corollary}[subsection]{Corollary}
\newtheorem{Theorem}[subsection]{Theorem}

\theoremstyle{definition}
\newtheorem{Definition}[subsection]{Definition}
\newtheorem{Remark}[subsection]{Remark}
\newtheorem{Example}[subsection]{Example}

\newcommand{\NEd}[1]{\begin{picture}(1.9,1)
\put(2.4,-1){\line(1,1){1}}
\put(3.6, .2){\vector(1,1){1}}
\put(2.7, 0.15){$\scriptstyle{#1}$}
\end{picture} &&} 

\newcommand{\NE}[1]{\begin{picture}(1.9,1)
\put(2.4,-1){\vector(1,1){2.2}}
\put(2.7, 0.15){$\scriptstyle{#1}$}
\end{picture} &&}

\newcommand{\SW}[1]{\begin{picture}(1.9,1)(-1.2,0.07)
\put(4.6 , 1.2){\vector(-1,-1){2.1}}
\put(2.7, 0.15){$\scriptstyle{#1}$}
\end{picture} &&} 

\newcommand{\SE}[1]{\begin{picture}(1.9,1)(0,0.0)
\put(2.4 , 1.2){\vector(1,-1){2.2}}
\put(1.9, -0.05){$\scriptstyle{#1}$}
\end{picture} &&}

\newcommand\bs{\backslash}

\newcommand\BC{\mathbb{C}}

\newcommand\BZ{\mathbb{Z}}
\newcommand\BR{\mathbb{R}} 

\newcommand\BS{\mathbb{S}}
\newcommand\BP{\mathbb{P}}

\def\thm{\begin{Theorem}}
\def\cor{\begin{Corollary}}
\def\mylemma{\begin{Lemma}}
\def\defn{\begin{Definition}}
\def\enddefn{\end{Definition}}
\def\epr{\end{Theorem}}
\def\ex{\begin{Example}}
\def\demo{\begin{proof}}
\def\enddemo{\end{proof}}
\def\rem{\begin{Remark}}
\def\Ref{\ref}
\def\goth{\mathfrak}

\DeclareMathOperator{\im}{Im}
 
\DeclareMathOperator{\h}{H} 
\DeclareMathOperator{\SU}{SU}
\DeclareMathOperator{\SO}{SO}
\DeclareMathOperator{\SL}{SL}
\DeclareMathOperator{\GL}{GL}
\DeclareMathOperator{\U}{U}
\DeclareMathOperator{\s}{S}

\def\by#1{{\bf #1}, }
\def\paper#1{{\em #1\/}, }
\def\yr#1{(#1)}
\def\book#1{{\em #1\/}, }
\def\pp#1{#1}
\def\publ{}
\def\publaddr{}
\def\vol{}
\def\jour{}
\def\Ref{\ref}
\def\issue#1{(#1)}

\begin{document}

\title{Nonzero degree tangential maps between\\dual symmetric spaces}
\asciititle{Nonzero degree tangential maps between dual symmetric spaces}
\author{Boris Okun}
\address{Department of Mathematical Sciences\\
	University of Wisconsin--Milwaukee\\
	Milwaukee, WI 53201, USA}  
\email{okun@uwm.edu}

\begin{abstract} 
We construct a tangential map from a locally symmetric space of
noncompact type to its dual compact type twin. By comparing the
induced map in cohomology to a map defined by Matsushima, we conclude
that in the equal rank case the map has a nonzero degree.
\end{abstract}

\primaryclass{53C35 }
\secondaryclass{57T15, 55R37, 57R99}

\keywords{Locally symmetric space, duality, degree, tangential map, 
Matsushima's map}

\maketitle

\section{Introduction}

It has been known since work of Matsushima \cite{M}, that there exists a map $j^*\co \h^*(X_u,\BR)\to \h^*(X,\BR)$, where $X=\Gamma\bs G/K$ is a compact locally symmetric space of noncompact type and $X_u$ is its global dual twin of compact type. Moreover, Matsushima showed that this map is monomorphic and, up to a certain dimension depending only on the Lie algebra of $G$, epimorphic. A refinement of Matsushima's argument, due to Garland \cite{G} and Borel \cite{B}, allowed the latter to extend these results to the case where $X$ is noncompact but only has a finite volume (but $\Gamma$ has to be an arithmetic group), thus giving, for example, a computation of the rational cohomology of $\SL(\BZ)\bs \SL(\BR)/\SO$, thus a computation of the rational algebraic K--theory of the integers. However, since the construction of the map $j^*$ is purely algebraic (in terms of invariant exterior differential forms), the natural question arises: is there a topological map $X \to X_u$ inducing $j^*$ in cohomology?  

For the case of $X$ being a complex hyperbolic manifold this kind of map was constructed by Farrell and Jones in \cite{FJ}, where it was used to produce nontrivial smooth structures on complex hyperbolic manifolds. In general, the answer to this question is negative, since Matsushima's map does not necessarily take rational cohomology classes into rational ones (this irrationality is measured by Borel regulators, see \cite{B1} for details), so it cannot be induced by a topological map. However, for the case when $G$ and $K$ are of the same rank we shall see that the answer is (virtually) positive.

In this paper we show that there is a finite sheeted cover $X'$ of $X$ and a tangential map (i.e.\ a map covered by a map of tangent bundles) $X' \to X_u$. In the case where $G$ and $K$ are of equal rank we prove that $k^*$ coincides with Matsushima's map  $j^*$ and therefore has nonzero degree. Thus, we are able to recover Farrell and Jones' result, but our method is different and more general.

We note that the construction of the map $k$ has some choices made in it.  Surprisingly enough this freedom does not appear in the equal rank case, but it is crucial for constructing a nonzero degree map in the general case.  It is proved in \cite{O} that the map $k$ can always be chosen to be of nonzero degree (but we cannot ensure that $k^*=j^*$). In a forthcoming paper, we show how this map can be used to produce new examples of exotic smooth structures on nonpositively curved symmetric spaces. 

The paper is organized as follows. In Section \ref {PRE}, we introduce some notation and give a quick overview of (the relevant part of) the theory of symmetric spaces. For reader's convenience, we describe four examples of dual symmetric spaces in Section \ref{EX}. In Section \ref {MM}, we recall Matsushima's construction. Section \ref{CTM} contains the construction of a tangential map. It turns out that in the equal rank case the induced map in cohomology is, in fact, Matsushima's map, thus it has a nonzero degree.  We prove this in Section \ref{ERC}.  

The material in this paper is part of the author's 1994 doctoral dissertation \cite{O} written at SUNY Binghamton under the direction of  F. T. Farrell. I am deeply grateful to him for bringing this problem to my attention and for invaluable advice during the course of this work. I also wish to thank S. C. Ferry for his constant encouraging attention and numerous useful discussions. The author was partially supported by NSF grant.

\section{\label{PRE} Notation and preliminary results}

We will use notation from \cite{B} and \cite{M}. Throughout this paper  $G$ will be a real semisimple algebraic linear Lie group (a subgroup of $\GL(n,\BR)$) and  $K$ will be its maximal compact subgroup. Let $G_c$ denote the complexification of $G$ and let $G_u$ denote a maximal compact subgroup of $G_c$.

\defn\label{DEF}  The factor spaces $G/K$ and $X_u=G_u/K$ are called dual symmetric spaces of noncompact type and compact type respectively. Let $\Gamma$ be a torsion-free discrete subgroup of $G$ of finite covolume. The space $X=\Gamma \bs G/K$ is a locally symmetric space of noncompact type. Slightly abusing terminology, we will also refer to  $X$ and $X_u$ as dual symmetric spaces.
\end{Definition}

We note that these spaces have natural Riemannian metrics, coming from the Killing form on Lie group $G$. It is a standard fact \cite{He} that for compact type, this metric has nonnegative curvature and for noncompact type, the curvature is nonpositive. In particular, it follows from Cartan--Hadamard theorem that the space $X$ is a $K(\Gamma ,1)$ manifold. Since $X$ has finite volume, it follows from the work in \cite{KM},\cite{BH},\cite{MT} and \cite{PR} that $X$ has finite homotopy type. (This also follows from a much more general theorem on p. 138 of \cite{BGS}.) 

The projections $p\co \Gamma\bs G\to  \Gamma\bs G/K$ and $p_u\co G_u \to G_u/K$ are principal fiber bundles with structure group $K$. They are called canonical principal bundles. These canonical bundles are classified by the maps to the classifying space $BK$ which we will denote by $c_p$ and $c_{p_u}$, respectively.

What we have said so far can be summarized in the following three commutative diagrams:
$$
\begin{CD}
K  @>\subset >>   G_u       \\
@VV\cap V             @VV\cap V    \\
G          @>\subset >>   G_c
\end{CD}
$$

\medskip

$$
\begin{CD}
K @= K && \hskip 12em && K @= K\\
@VVV @VVV &&  @VVV @VVV \\
\Gamma \bs G  @>>> EK && && G_u @>>> EK \\
@V{p}VV @VVV &&  @Vp_uVV @VVV \\
\hskip -3em X=\Gamma \bs G/K @>c_p>> BK &&  && \hskip -3em X_u=G_u/K @>c_{p_u}>> BK
\end{CD}
$$

 The Lie algebras of $G$ and $G_u$ have dual Cartan decompositions
\begin{equation}\label{cartan}
\begin{split}
\goth{g}=\goth{k}\oplus\goth{p}\\  
\goth{g}_u=\goth{k}\oplus i\goth{p}
\end{split}
\end{equation}
where $\goth{k}$ denotes the Lie algebra of the group $K$.

It is well known that tangent bundles of symmetric spaces can be expressed in terms of Lie algebras and canonical bundles. We have the following Lemma:
 
\mylemma[see, e.g., pp. 209--210 of \cite{GHV})\label{t}]The tangent bundles of the dual symmetric spaces have the following  expressions: 
$$\tau(X)=p\times_K\goth{p}$$ 
$$\tau(X_u)=p_u\times_K i\goth{p}$$ 
Here we consider $\goth{p}$ and $i\goth{p}$ as $K$--modules via adjoint action of $K$.
\end{Lemma}
 
Lemma \ref{t} implies: 

\mylemma{\label{Tang}} Any map $g\co X \to X_u$ between the dual symmetric spaces which preserves canonical $K$--bundle structure is tangential.
\end{Lemma}

\demo
The proof is obvious, since $\goth{p}$ and $i\goth{p}$ are isomorphic as 
$K$--modules.
\enddemo

\section{\label{EX} Examples}

In this section, we give some examples of dual symmetric spaces, which often appear in geometric topology.

\ex{} Let $G=\SL(n,\BR)$. Then the maximal compact subgroup $K$ of $G$ is $\SO(n,\BR)$. The complexification $G_c$ of $G$ is $SL(n,\BC)$ and its maximal compact subgroup is $SU(n)$, so we have
$$X_u=\SU(n)/\SO(n,\BR).$$
As a discrete subgroup $\Gamma \subset G$ one would like to take the group $\SL(n,\BZ)$, unfortunately this group has torsion, so it does not quite fit in our setup. We get around this problem by recalling a result of Selberg that there is a torsion-free subgroup
$\Gamma \subset \SL(n,\BZ)$ which has a finite index. Thus, we have 
$$X=\Gamma \bs \SL(n,\BR) / \SO(n,\BR).$$
We note that primitive elements in cohomology of $X$ (after passing to  the limit as $n\to \infty$) give, at least rationally, the algebraic K--theory of integers.
\end{Example}

\ex\label{RH} This example shows the duality between hyperbolic manifolds and spheres.
 
Let  $q(x)$ denote the nondegenerate indefinite quadratic form on $\BR^{n+1}$ defined by
$$ q(x)=x_1^2+x_2^2+\dots +x_n^2-x_{n+1}^2 $$
where $x \in \BR^{n+1}$ and $x_i$ denotes the $i$-th coordinate of $x$.

Let $G=\SO^+(n,1,\BR)$ denote the connected component of the identity of the Lie group of all isometries of the form $q$ which have determinant equal to 1.  The maximal compact subgroup $K$ of $G$ in this case is $\SO(n,\BR)$. The manifold $G/K=\SO^+(n,1,\BR)/\SO(n,\BR)$ is a (real) $n$--dimensional hyperbolic space. If $\Gamma$ is a torsion free discrete subgroup of $G=\SO^+(n,1,\BR)$, then the locally symmetric space of noncompact type $X=\Gamma\bs G/K=\Gamma\bs \SO^+(n,1,\BR)/\SO(n,\BR)$ is a (real) hyperbolic manifold.

Now let us look at the dual twin of $X$. The complexification $G_c$ of $G$ is $\SO(n,1,\BC)$. However, since all quadratic forms are equivalent over the field of complex numbers, we have $G_c=\SO(n+1,\BC)$. The maximal compact subgroup $G_u$ of $G_c$ is $\SO(n+1,\BR)$, so we have 
$$ X_u=\SO(n+1,\BR)/\SO(n,\BR)=\BS^{n}. $$
\end{Example}

\ex\label{CH}
This example is similar to the previous one, but instead of using a quadratic form over reals, we use a Hermitian form over complex numbers. This leads to duality between complex hyperbolic spaces and complex projective spaces.

Let $b(,)$ be  the nondegenerate indefinite Hermitian form on $\BC^{n+1}$ defined by
$$ b(x,y)=x_1\bar y_1+x_2\bar y_2+\dots +x_n\bar y_n-x_{n+1}\bar y_{n+1} $$
where $x,y \in \BC^{n+1}$ and $x_i$ denotes the $i$-th coordinate of $x$.

Let $G=\SU(n,1)$ denotes the Lie group of all isometries of the form $b(,)$ which have determinant equal to 1. The maximal compact subgroup $K$ of $G$ in this case is $\s(\U(n)\times \U(1))$. The manifold $G/K=\SU(n,1)/\s(\U(n)\times \U(1))$ is a complex $n$--dimensional hyperbolic space.  If $\Gamma$ is a torsion-free discrete subgroup of $G=\SU(n,1)$ then locally symmetric space of noncompact type $X=\Gamma\bs G/K=\Gamma\bs \SU(n,1)/\s(\U(n)\times \U(1))$ is a complex  hyperbolic manifold.

On the compact side we have
$$G_c=\SL(n+1,\BC)$$
$$G_u=\SU(n+1)$$
$$X_u=\SU(n+1)/\s(\U(n)\times \U(1))=\BC\BP^n.$$
\end{Example}

\begin{Example}\label{Complex} Suppose $G=G_c$ is a complex semisimple Lie group, and let $G_u$ denote its maximal compact subgroup. The complexification of $G_c$ is then isomorphic to $G_c\times G_c$:
$$(G_c)_c=G_c\times G_c$$
The maximal compact subgroup of the complexification is $G_u\times G_u$, so we have dual symmetric spaces:
$$X=\Gamma\bs G_c/G_u$$
$$X_u=G_u=G_u\times G_u/G_u$$
In particular, if $G=\SL(n,\BC)$, we have a pair:
$$X=\Gamma\bs \SL(n,\BC)/\SU(n)$$
$$X_u=\SU(n)$$
\end{Example}

\section{\label{MM} Matsushima's map}

In this section we recall Matsushima's construction \cite{M} of a map in real cohomology $j^*\co \h^*(X_u,\BR) \to \h^*(X,\BR)$.

Matsushima's construction uses differential forms. If a group $H$ acts on a smooth manifold $Y$ we let $\Omega^H(Y)$ denote the complex of $H$--invariant differential forms on $Y$. 

 Let $X=\Gamma\bs G/K$ and $X_u=G_u/K$ be dual symmetric spaces. We have natural left actions of $G$ on $G/K$ and $G_u$ on $X_u$. It is well-known that $\Omega^G(G/K)$ and $\Omega^{G_u}(X_u)$ consist of harmonic forms, and we have
\begin{equation}\label{1}
\begin{split} 
\Omega^G(G/K)=\h^*(\goth{g},\goth{k}) \\
 \Omega^{G_u}(X_u)=\h^*(\goth{g}_u,\goth{k})
\end{split}
\end{equation}
where the right-hand side terms denote relative Lie algebra cohomology.

Using Cartan decomposition \eqref{cartan} we see that 
\begin{equation}\label{2}
\h^*(\goth{g},\goth{k})=\h^*(\goth{g}_u,\goth{k}).
\end{equation}
Since $X_u$ is compact closed manifold it follows from Hodge theory that 
\begin{equation}\label{3} 
\Omega^{G_u}(X_u)=\h^*(X_u,\BR).
\end{equation}

Combining equations \eqref{1}, \eqref{2} and \eqref{3} we get 
$$\Omega^G(G/K)= \h^*(X_u,\BR).$$
Now let us consider the complex $\Omega^\Gamma(G/K)$  of $\Gamma$--invariant forms. The projection $G/K \to \Gamma\bs G/K$ induces the isomorphism
$\Omega^\Gamma(G/K)=\Omega(X)$ and therefore we have 
$$ \h^*(X,\BR)= \h^*(\Omega^\Gamma(G/K)).$$
Since the elements of $\Omega^G(G/K)$ are closed forms, the inclusion  $\Omega^G(G/K) \subset \Omega^\Gamma(G/K)$ induces a homomorphism 
$$j^*\co \h^*(X_u,\BR) \to \h^*(X,\BR),$$
which we will be calling Matsushima's map.

\thm{\rm \cite{M}\label{MT}}\qua Let $X=\Gamma \bs G/K$ and $X_u=G_u/K$ be dual symmetric spaces. Assume $X$ to be compact, then
\begin{enumerate}
\item $j^*$ is injective.

\item $j^q$ is surjective for $q\leq m(\goth{g})$, where $m(\goth{g})$ is some constant, explicitly defined in terms of Lie algebra $\goth{g}$ of $G$.
\end{enumerate}
\end{Theorem}

\rem
Since $X$ is compact, the injectivity part follows from Hodge theory. In \cite{M} Matsushima is mostly concerned with the surjectivity of $j^*$.
\end{Remark}

Matsushima's map is closely related to the characteristic classes, as can be seen from the following lemma.

\mylemma[\cite{B1},\cite{KT}\label{tria}] Let $X=\Gamma \bs G/K$ and $X_u=G_u/K$ be dual symmetric spaces. Then the  diagram
$$ 
\begin{CD} 
&& \h^*(X_u,\BR) \\
\SW{j^*} @AA{{c_{p_u}}}^*A \\
\h^*(X,\BR) @<<{{c_p}^*}< \h^*(BK,\BR) 
\end{CD}
$$
commutes. 
\end{Lemma}
 
\section{\label{CTM} Construction of tangential map}

In this section, we show the existence of a tangential map between dual symmetric spaces. 

\thm{\label{ta}} Let $X=\Gamma\bs G/K$ and $X_u=G_u/K$ be dual symmetric spaces.
Then there exist a finite sheeted cover $X'$ of $X$ (i.e.\ a subgroup
$\Gamma'$ of finite index in $\Gamma$, $X'=\Gamma '\bs G/K$) and a tangential map $k\co X' \to X_u$. 
\end{Theorem}

\demo Consider  the canonical principal fiber bundle with structure group $K$ over $X$ $p\co \Gamma\bs G\to  \Gamma\bs G/K$. If we extend the structure group to the group $G$ we get a flat principal bundle: $\Gamma\bs G\times_K G=G/K\times_\Gamma G\to  \Gamma\bs G/K$ \cite{KT}. Extend the structure group further to the group $G_c$. The resulting bundle is a flat bundle with an algebraic linear complex Lie structure group, so by a theorem of Deligne and Sullivan \cite{DS} there is a finite sheeted cover $X'$ of $X$ such that the pullback of this bundle to $X'$ is trivial. This means that for $X'$ the bundle obtained by extending the structure group from $K$ to $G_u$ is trivial too, since $G_u$ is the maximal compact subgroup of $G_c$. Consider now the following diagram:
$$
\begin{CD} 
&& X_u \\
\NEd{k} @VV{c_{p_u}}V \\
X' @>{c_{p'}}>> BK \\
\SE{\simeq 0} @VViV \\ 
&& BG_u
\end{CD} 
$$ 
Here the map $c_{p'}$ is a classifying map for the canonical bundle $p'$ over $X'$.

The map $i$, induced by standard inclusion $K\subset G$ is a fibration with a fiber $X_u=G_u/K$. Note that the inclusion of $X_u$ in $BK$ as a fiber also classifies canonical principal bundle $p_u\co G_u\to G_u/K$. By the argument above the composition $i\circ c_{p'}\co X'\to BG_u$ is homotopically trivial. 

Choose a homotopy contracting this composition to a point. As the map $i$ is a fibration we can lift this homotopy to $BK$.  The image of the end map of the lifted homotopy is contained in the fiber $X_u$, since its projection
to $BG_u$ is a point. Thus we obtain a map $k\co X'\to X_u$ which makes the upper triangle of the diagram homotopy commutative. It follows that the map $k$ preserves canonical bundles on the spaces $X'$ and $X_u$:
$k^*(p_u)=p'$. By Lemma \Ref{Tang}, the map $k$ is tangential.
\enddemo

\section{\label{ERC} Nonzero degree in equal rank case}

In this Section we prove that if Lie groups $G_u$ and $K$ have equal ranks then the map constructed in Theorem \Ref{ta} induces Matsushima's map in cohomology and therefore has a nonzero degree in compact case.

\defn The dimension of the maximal torus of a compact Lie group $H$ is called the rank of $H$.
\end{Definition}

The following lemma is well known.

\mylemma\label{Eqr}\qua{\rm \cite{Bo0}}\qua
 Let $H$ be a compact Lie group and $K$ be its Lie subgroup. If the ranks of $H$ and $K$ are equal then the map classifying canonical bundle on the homogeneous space $H/K$ induces epimorphism in rational cohomology. 
\end{Lemma}

\mylemma{\label{Agree}} Let   $g\co X \to X_u$ be any between the dual symmetric spaces which preserves canonical $K$--bundle structure. Then 
$$g^*|_{\im c^*_{p_u}}=j^*|_{\im c^*_{p_u}}.$$
\end{Lemma}

\demo The condition implies that triangle
$$
\begin{CD}
&& X_u \\
\NE{g} @VV{c_{p_u}}V \\
X @>>{c_p}> BK \\
\end{CD}
$$
is homotopy commutative. Passing to cohomology we see that the diagram 
$$ 
\begin{CD} 
&& \h^*(X_u,\BR) \\
\SW{g^*} @AA{{c_{p_u}}}^*A \\
\h^*(X,\BR) @<<{{c_p}^*}< \h^*(BK,\BR) 
\end{CD}
$$
commutes. Comparing this diagram to Lemma \Ref{tria} we see that $g^*$ and $j^*$  have to coincide on the image of the map $c_{p_u}^*$. 
\enddemo

\thm{\label{tb}}
Let $X=\Gamma\bs G/K$ and $X_u=G_u/K$ be dual symmetric spaces.
Let $X'$ be the finite sheeted cover of $X$ and $k\co X'\to X_u$ be the tangential map, constructed in Theorem \Ref{ta}.
If the groups $G_u$ and $K$ are of equal rank then the map induced by $k$ in cohomology coincides with Matsushima's map $j^*$. 
\end{Theorem}

\demo
Since $k$ preserves canonical bundles, by Lemma \Ref{Agree} the two maps agree on the image of the map $c_{p_u}^*$. By Lemma \Ref{Eqr} the map $c_{p_u}^*$ is epimorphic so $k^*=j^*$.
\enddemo

\cor
If the groups $G_u$ and $K$ are of equal rank and the group $\Gamma$ is cocompact in $G$ then the map $k$ has a nonzero degree.
\end{Corollary}

\demo
This follows from Theorem \Ref{MT} and Theorem \ref{tb}.
\enddemo

Among examples considered above, the examples \ref{RH} for even $n$, and \ref{CH} for all $n$ are of equal rank.

\Addresses\recd

\begin{thebibliography}{10}

\bibitem{BGS}
\by{W Ballmann, M Gromov\ and\ V Schroeder}
 \book{Manifolds of nonpositive curvature} 
\jour{Progr. Math.}
\vol{61}
\publ{Birkh\"auser}
\publaddr{Boston}
\yr{1985}


 \bibitem{Bo0}
 \by{A Borel}
\paper{Sur la cohomologie des espaces fibr\'{e}s principaux et des espaces
homog\'{e}nes de groupes de Lie compacts}
\jour{Ann. of Math.}
\vol{57}
\issue {2}
\yr{1953}
\pp{115--207}


  \bibitem{B}
\by{A Borel}
\paper{Stable real cohomology of arithmetic groups}
\jour{Ann. Sci. Ec. Norm. Sup.}
\vol{7}
\issue {4}
\yr{1974}
\pp{235--272}


 \bibitem{B1}
 \by{A Borel}
\paper{Cohomologie de $\SL_n$ et valeurs de fonctions zeta aux points entiers}
\jour{Ann. Sc. Norm. Super. Pisa, Cl. Sci.}
\vol{4}
\issue{4}
\yr{1977}
\pp{613--636}

 \bibitem{BH}
 \by{A Borel and Harish-Chandra}
\paper{Arithmetic subgroups of algebraic groups}
\jour{Ann. of Math.}
\vol{76}
\yr{1962}
\pp{485--535} 

  \bibitem{DS}
 \by{P Deligne and D Sullivan}
\paper{Fibr\'{e}s vectoriels complexes \`{a} groupe structural
 discrete}
\jour{C. R. Acad. Sc. Paris, S\'{e}r. A}
\vol{281}
\yr{1975}
\pp{1081--1083}

  \bibitem{FJ}
 \by{F\,T Farrell and L\,E Jones}
\paper{Complex hyperbolic manifolds and exotic smooth structures}
\jour{Invent. Math.}
\vol{117}
\yr{1994}
\pp{57--74}

  \bibitem{G}
 \by{H Garland}
\paper{A finiteness theorem for $K_{2}$ of a number field}
\jour{Ann. of Math.}
\vol{94}
\issue {2}
\yr{1971}
\pp{534--548}

  \bibitem{GHV}
 \by{W Greub, S Halperin and R Vanstone}
\book{Connections, curvature, and cohomology}
\vol{II}
\publ{Academic Press}
\publaddr{New York and London}
\yr{1973}
 
  \bibitem{He}
 \by{S Helgason}
\book{Differential geometry and symmetric spaces}
\publ{Academic Press}
\publaddr{New York and London}
\yr{1962}

  \bibitem{KT}
 \by{F Kamber and P Tondeur}
\book{Foliated bundles and characteristic classes}
\publ{Springer-Verlag}
\publaddr{Berlin-New York}
\yr{1975}

 \bibitem{KM}
 \by{D\,A Kazhdan and G\,A Margulis}
\paper{A proof of Selberg's hypothesis}
\jour{Math. Sbornik}
\vol{75}
\issue{117}
\yr{1968}
\pp{162-168}

  \bibitem{M}
 \by{Y Matsushima}
\paper{On Betti numbers of compact, locally symmetric Riemannian
 manifolds}
\jour{Osaka Math. J.}
\vol{14}
\yr{1962}
\pp{1--20}

 \bibitem{MT}
 \by{G\,D Mostow\ and\ T Tamagawa}
\paper{Arithmetic subgroups of algebraic groups}
\jour{Ann. of Math.}
\vol{76}
\yr{1962}
\pp{446--463} 

  \bibitem{O}
\by{B Okun}
\paper{Nonzero degree maps between dual symmetric spaces}
\publ{PhD Thesis}
\publaddr {SUNY Binghamton}
\yr{1994}  

 \bibitem{PR}
 \by{G Prasad\ and\ M\,S Raghunathan}
\paper{Cartan subgroups and lattices in semi-simple groups}
\jour{Ann. of Math.}
\vol{96}
\yr{1972}
\pp{296--317} 

\end{thebibliography}
\end{document}